\documentclass[11pt]{amsart}
\usepackage{amsfonts,amsmath,amssymb,eucal,tipa,stmaryrd}
\usepackage[all]{xy}
\usepackage{textcomp}
\usepackage{hyperref}
\usepackage[utf8]{inputenc}
\usepackage{shuffle}
\usepackage[T1]{fontenc}
\usepackage[czech,polish]{babel}

\begin{document}

\parindent=0pt
\parskip=6pt

\newcommand{\ie}{\textit{ie}\,}
\newcommand{\eg}{\textit{eg}\,}
\newcommand{\se}{{\sf e}}
\newcommand{\cf}{{\textit{cf}\,}}
\newcommand{\nsymm}{{\sf NSymm}}
\newcommand{\qsymm}{{\sf QSymm}}
\newcommand{\symm}{{\sf Symm}}
\newcommand{\mH}{{\mathfrak H}}
\newcommand{\vc}{{\copyright}}
\newcommand{\zz}{{\zeta^{\mathfrak m}}}
\newcommand{\Gl}{{\rm Gl}}
\newcommand{\oh}{{\mathfrak o}}
% {{\textcent}}
\newcommand{\Br}{{\mathbb B}}
\newcommand{\C}{{\mathbb C}}
\newcommand{\D}{{\rm Diff}}
\newcommand{\h}{{\mathfrak h}}
\newcommand{\PGl}{{\rm PGl}}
\newcommand{\Pic}{{\rm Pic}}
\newcommand{\Q}{{\mathbb Q}}
\newcommand{\R}{{\mathbb R}}
\newcommand{\Sl}{{\rm Sl}}
\newcommand{\Sf}{{\rm SF}}
\newcommand{\SU}{{\rm SU}}
\newcommand{\T}{{\mathbb T}}
\newcommand{\U}{{\mathbb U}}
\newcommand{\wzw}{{\bf wzw}}
\newcommand{\Z}{{\mathbb Z}}
\newcommand{\G}{{\mathbb G}}

\newcommand{\Cz}{{\sf Cos}}
\newcommand{\CK}{{\rm CK}}
\newcommand{\CM}{{\rm CM}}
\newcommand{\Lc}{{\mathcal L}}
\newcommand{\g}{{\mathfrak g}}
\newcommand{\bl}{{\boldsymbol{\lambda}}}
\newcommand{\cz}{{\check{c}}}
\newcommand{\bro}{{\boldsymbol{\rho}}}
\newcommand{\st}{{\sf st}}
\newcommand{\bvdash}{{\boldsymbol{\vdash}}}
\newcommand{\om}{{\Omega S^2}}

\newcommand{\m}{{\rm max}}
\newcommand{\dvol}{{d{\rm vol}}}
\newcommand{\dR}{{\rm dR}}
\newcommand{\Maps}{{\rm maps}}
\newcommand{\haut}{{\mathcal{H}{\it aut}}}
\newcommand{\cD}{{\mathcal D}}
\newcommand{\HL}{{{\mathcal H}L}}
\newcommand{\Top}{{\rm Top}}

\title{A homotopy-theoretic context for CKM/Birkhoff renormalization}

\author[Jack Morava]{Jack Morava}

\address{Department of Mathematics, The Johns Hopkins University,
Baltimore, Maryland}

\begin{abstract}{We propose a geometric object slightly subtler than a complex line bundle with connection of some kind, a two-sphere fibration with structure group $\Omega^2_e S^2$ parametrizing the space of dimensional regularizations $\Z \ni d \mapsto d - s \in {\mathbb P}_1(\C)$, as in \cite{1, 12, 13} and throughout the general metaphysics of renormalization theory.}\end{abstract}

\maketitle 

{\bf \S I}\begin{footnote}{Tangential comments are collected in an appendix. In this section (co)homology has coefficients $\Z$, unless not.}\end{footnote}  

% dated 29Sept023

The object of this project is to construct an arithmetic de Rham cohomology \cite{15,19,33,38} for (some parts of) quantum field theory, interpreted geometrically along Chern-Weil lines, with coefficients defined over a suitable ring of multi-zeta values: roughly, a theory of renormalized geometric cycles \cite{36}(\S 8) in classical quantum electrodynamics in which the Feynman integrals can be coherently evaluated in terms of multizeta values. [This is not true for many standard physical models.] \bigskip

{\bf 1.0} The homomorphism 
\[
\bl : \T \to \D(\R^2_+; 0,\infty) ( \subset \Omega^2_e S^2) 
\]
of topological groups defined by plane rotation ($\T \cong {\rm SO}(2)$) is homotopically a maximal torus; it induces a map 
\[
B\bl : B\T \to \Omega S^2 
\]
of spaces (though not of $H$-spaces). Its domain is the group completion of the free commutative monoid of points on the projective line, while the target is the group completion of Segal's monoid of its virtually finite subsets\begin{footnote}{$V \subset \C P_1 = H^+ \cup H^- \Rightarrow \#(V \Delta H^\pm) < \infty$}\end{footnote}, as in the Goto-Nambu model. 

Because $B\T$ is simply-connected there is a lift 
\[
\xymatrix{
{} & \Omega S^3 \ar[d]^\Z \ar[r]^\cong & \Omega \SU(2) \ar[d]^\cong \ar[r] & \PGl_\C({\mathbb H}) \ar[d] \ar[r]^\simeq & 
B\T \\
B\T \ar[ur] \ar[r] & \Omega S^2 \ar@{.>}[r]^\nexists & J(S^2) = (\Omega \Sigma) (S^2) \ar[r] &  \vee_{n \geq 0} S^{2n} \ar[r] & S^2 \ar[u] }
\]
where IM James' construction
\[
X \to J(X) := (\Omega \Sigma) (X) \sim \vee_{n \geq 0} X^n
\]
(up to suspension) embeds a space in the free topological monoid it generates; see \cite{25}. The nonexisting dotted map tries to extend an $E_1$ to an $E_2$ structure; in this picture instead of thinking of divisors as sums of points with multiplicity (zeroes and poles) we imagine them (in low dimensions) as formal products of immersed two-spheres.\bigskip

{\bf 1.1 Definition}  A \l{}ine bundle on a based space $X$ is a map
\[ 
\lambda : X \to \Omega S^2 
\]
classifying a principal $\Omega^2_e S^2 := \wzw$ bundle, with fiber interpretable as the configuration space $\Omega^2 \R^3_+$ of the  simplest Wess-Zumino-Witten model \cite{36}(\S 6.3)). There is a Picard group \cite{50}
\[
\Pic_\wzw (X) := \pi_0 {\rm Maps}(\Sigma X;S^2) := \pi^2(\Sigma X) \to \Pic_\C(X)
\]
of equivalence classes of such things.

Lifting to universal covers defines a complex line bundle $\tilde{\lambda}: \tilde{X} \to \tilde{\Omega} S^2 \to \PGl_\C({\mathbb H}) \simeq B\T$ over a virtual cyclic cover of $X$ \cite{5}(\S 5.2). In \S III {\it iii}) below we propose to interpret $\lambda$ as a kind of homotopy-theoretic Cartan connection on a (dimensional regularization) two-sphere bundle over $X$; see \cite{12}(Th 5.1), \cite{13, 20, 42}, and further under digression {\it ii}). \bigskip

{\bf 1.2} Because $\Omega S^2$ is an $H$-space, iterating $B\lambda$ defines an action 
\[
J(B\T) \wedge \Omega S^2 \to \Omega S^2
\]
\ie 
\[
\xymatrix{
{} & J(B\T) \ar[d] \\
B\T \ar[r] \ar@{.>}[ur]  & \Omega S^2 }
\]
of the free topological monoid on $B\T$ upon $\Omega S^2$, making it something like a determinant line. On cohomology there is an induced action 
\[
H_*J(B\T) \otimes H^*\Omega S^2 \to H^*\Omega S^2
\]
of the noncommutative symmetric functions $\nsymm$  \cite{10}(\S 4.1F). Equivalently, the coaction
\[
H^*(\Omega S^2) \to H^*(\Omega S^2) \otimes H^*J(B\T) 
\]
makes $H^*(\Omega S^2)$ a comodule over Baker and Richter's dual Hopf algebra \cite{3} $\qsymm$ of quasi-symmetric functions (under the quasi-shuffle $*$ product induced by their representation $J(B\T) \to \U$) as below.

The induced  morphism
\[
H^*B\lambda : H^*(\Omega S^2)  \cong E(e) \otimes \Z[\gamma^*b] \to \Z[c] \cong H^*B\T
\]
of algebras sends the Bockstein element $e$ to 0 and $\gamma^k b \to c^k$ : Equivalently, writing $\gamma^k \bar{c} \in H_{2k} B\T$ for the dual to $c^k$, the action 
\[
H_*B\T \otimes H^*(\Omega S^2) \to H^*(\Omega S^2) ,
\]
can be written $\gamma^k \bar{c} \otimes \gamma^l b \to \gamma^{l-k} b$. We extend this convention below to coefficients in $M\U \otimes \Q$.\bigskip

 {\bf Definition} \cite{9}(\S 2.7 eq 125) :
\[
\cz (b) := \sum_{k \geq 0} \gamma^k b \cdot T^k  = \exp (bT) \in (H^*\Omega S^2)[[T]]
\]
is an analog of a generalized Chern character for such \l{}ine bundles. \bigskip

{\bf \S \; II over} $\Q$\bigskip

{\bf 2.1.0} The Boardman-Hurewicz ring homomorphism injects the Thom-Milnor ring $M\U_* \to H_*B\U$ of complex-oriented manifolds in the Pontryagin ring of $B\U$ and thus in the ring of symmetric functions, sending Chern-Weil integrals $c_I$ indexed by elementary symmetric functions $e_I$ to themselves \cite{35}(\S 1.2) indexed by the monomial symmetric functions $m_I$. The Hurewicz image of the universal formal group is defined over $\Z$ by 
\[
X,Y \mapsto m^{-1}(m(x) + m(Y)) \in \symm_*[[X,Y]]
\]
(with $m(X) = \sum_{i \geq 0} m_i x^{i+1}$ \cite{24}). This defines a canonical Chern-Dold character
\[
ch_{MU} : M\U^*(X) \otimes \Q \to H^*(X; \symm_* \otimes \Q) 
\]
which we will take to be an identification. The proalgebraic monoid $\Gamma$ of formal series $t(T) = \sum_{i \geq 0} t_i T^{i+1}$ under composition acts by
\[
(t \circ e_*) (T)  = t \circ \prod (1 +  x_i T) \mapsto \prod (1 + x_i t(T)) \in \symm_*[[T]]
\]
on the symmetric functions, defining the Landweber-Novikov operations in cobordism theory. \bigskip

{\bf 2.1.1} In Michael Hoffman's free graded algebra $\h = \Q\langle x,y \rangle$ of noncommutative polynomials in two variables, the sub-vector space 
\[
\h^1 := 1 + \h \cdot y
\]
generated by monomials ending in $y$ is in fact a subalgebra: which can be identified with the free graded algebra 
\[
\nsymm := \Q \langle y_1,\dots \rangle, y_k= x^{k-1}y
\]
of noncommutative symmetric functions \cite{10}(\S 4.1F). 

As a vector space, $\h^1$ has a remarkable family of commutative associative quasi-shuffle algebra products [$*, \star, \shuffle, \dots$, \cite{15}(\S 2)] on $\qsymm$, all of them polynomial \cite{9}(\S 2.7), \cite{8}(\S 1.7.4) over $\Q$ on generators indexed by the Lyndon words in two letters \cite{44}. \bigskip

{\bf 2.1.2} Baker and Richter's multiplicative cohomology theory $M\xi$ is defined by the Thom ring-spectrum associated to their morphism
\[
\xi \sim J(B\T) \to \U .
\]
In this language \cite{9}(\S 3.4  eq 155), they show :

{\bf Theorem} {\it There is a canonical natural equivalence 
\[
\mathfrak{br}_* : M\xi^*(X) \otimes \Q \to H^*(X;\h^1_*) 
\]
of multiplicative cohomology theories for finite CW spaces, with the subscript indicating the $*$ product on $\h^1$, (graded Hopf) dual to the composition product on $\nsymm$.} \bigskip

It is useful to define a regularization isomorphism 
\[
\h^1 \cong \h^0[T], \h^0 = 1 + x \cdot \h \cdot y 
\]
where $y \to y_1 = x^0 y \to T$, but $T$ is central -- whereas $y_1$ isn't \cite{9}(\S 2.4 eq 96). 

More generally, in nonsimply-connected contexts \cite{5}, \cite{8}(\S 3.10.2, def 3.402 fig 22),\cite{45} \dots one encounters Alexander's algebra 
\[
\Lambda := \Z[\pi_1(\T)] = \Z[T^{\pm 1}]
\]
suggesting cohomology with twisted coefficients in a local system \cite{22}{\S 5.2)
\[
X \mapsto H^*({\widetilde X},\Lambda) \otimes \h^1
\]
\ie  with
\[
\bvdash : T \to T - \gamma \ni \R
\]

\cite{9}(\S 2.8 eq 122) as an interesting functor. In what follows we assume for simplicity \cite{31}(\S 2.3.1) that 
$\gamma/2\pi$ is rational.\bigskip

{\bf Proposition} {\it The commutative diagram 
\[
\xymatrix{
M\xi^* \otimes \Q \ar[d] \ar[r] & \h^0_*[T] \\
MU^*\otimes \Q \ar[r]^\cong & \symm_* \otimes \Q \ar[u]_\cong}
\]
identifies the image of the abelianization ring homomorphism  $M\xi \otimes \Q \to M\U \otimes \Q$ with inclusion $\symm \subset \h^1_*$ of the symmetric functions as infinitesimal $\diamondsuit$-algebra} \cite{18}(Th 4.30, \cite{10}(\S 4.1G).\bigskip

This makes $\h^0[T] \cong M\xi \otimes \Q$ an $M\U \otimes \Q$ - algebra, even though $M\xi$ does not seem to be an $M\U$ -module spectrum.\bigskip

{\bf \S 2.3 Definition} Cartier's dimensional renormalization flow has operator generator
\[
 \log_{M\U} (\partial_T) = \sum_{n \geq 1} \frac{\C P_{n-1}}{n} \cdot \partial_T^n \in {\rm End}(\h^0_*[T])
\]
\cite{9}(\S 2.7 eq 129), with $[\C P^k]$ mapping to the symmetric power $p_k$ by MacDonald's $L$ \cite{23}(\S I.6 Ex 10). \bigskip

In complex cobordism, we have a coaction 
\[
M\U^*\Omega S^2 \to  M\U^* \Omega S^2 \otimes_{MU} M\U^* J(B\T) \to M\U^*\Omega S^2 \otimes_\symm \h^1_*
\]
induced by the action of $J(B\T)$ on $\Omega S^2$.\bigskip

{\bf Proposition} : {\it The generalized Chern class 
\[
\cz_{\rm MU}(b) := \exp_{*/{\rm MU}}(bT)
\] 
maps under renormalization \cite{9}(\S 2.7), \cite{8}(\S 1.7.8) to the Schneider-Teitelbaum modulus $\st_{MU}(b)$ \cite{39}
\[
Z_\shuffle \cz_{MU}(b) = \exp_*(- \log_{MU} (\partial_T)) \cdot \exp_* (bT) =  \exp_*(- \log_{MU}(b)).
\]
of the Quillen-Lazard formal group law. }\bigskip

Here renormalization $Z_\shuffle$ is defined by a change of quasi-shuffle product, with $M\U \otimes \Q$ as infinitesimal $\diamondsuit$-algebra (which is preserved by the Cartier flow). The argument follows \cite{9}(\S 2.4 eq 105), \ie d'apres Sir Isaac,
\[
\sum \lambda_r t^r = \exp \sum (-1)^{r-1}\psi_r \frac{t^r}{r} .
\] 

{\bf NB} Quasi-triangular quasi-Hopf algebra structures on 
\[
M\U^* J(B\T) \cong \symm^* \otimes \qsymm^*
\]

look to be very rich: \cite{9}(\S 3.4 eq 162), \cite{8}(Prop 3.393,\S 3.10.6) $\Box$\.] \bigskip

{\bf \S III  CKM renormalization in QFT, eg QED}\bigskip

Birkhoff renormalization, following Connes, Kreimer, Marcolli \cite{12, 13} \dots, defines a composition 
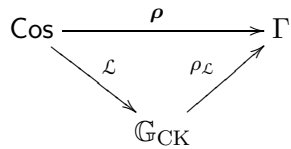
\[
\xymatrix{
\Cz \ar[dr]^\Lc \ar[rr]^\bro & {} & \Gamma \\
{} & \G_\CK \ar[ur]^{\rho_\Lc} }
\]
of group homomorphisms, where:  \bigskip

{\it i}) $\Cz$ is a group object classifying a cosmic Tannakian category of vector bundles over the punctured disk, endowed with flat equisingular connections \cite{13}(Th 2.25, esp lemma 7): very roughly, Riemann-Hilbert problems (\eg Yukawa) over $\C$ defined by connections
\[
\nabla \sim \partial_z + \sum \frac{A_i}{z - z_i} .
\]
As a topological group it is infinite-dimensional pro-unipotent with graded Lie algebra free on generators in consecutive positive degrees. Its exponential map is invertible and it is convenient to think of the free pro-unipotent associative algebra of functions on $\Cz$ as the universal enveloping algebra of its Lie algebra, cf remarks below. The continuous cochain Lie cohomology of such objects is a reasonable candidate \cite{8}(\S A.280) for $H^*B\Cz \otimes \Q$. \bigskip

{\it ii}) $\Lc$ signifies the choice of a Lagrangian action functional for a classical quantum field theory such as electrodynamics or the Salam-Weinberg \SU(2) model \cite{46} of weak interactions. 

The group $\G_\CK$ is similarly infinite-dimensional, with a locally finitely generated Lie algebra of diffeographisms defined by dismembering $\exp(-\Lc)$ into something like an asymptotic expansion with terms indexed by weighted (Feynman) graphs, using Wick's theorem \cite{20}. There is a (pre-) sheaf of QFT models \cite{1},\cite{13}(\S 2.2.3) over the category of Riemannian $d$ - manifolds $(X,g)$ further endowed with meromorphic families $\Delta^{-s}$ ((with fixed simple poles on the real axis \cite{27,42}) of Schwarz kernels for the conformal Laplacian \cite{30}(\S 1.3)) on a two-sphere bundle over $X$: with fibers parameterizing deformations of the purportedly physical dimension $d$. \bigskip

{\it iii}) Birkhoff renormalization \cite{10}(\S 5) associates to $\G_\CK$, a renormalization group flow defined by the choice of a smooth Jordan curve encircling the origin in $\C P_1$, an element of the group \cite{28}(\S 1.0.2) of formal automorphisms of the line at the origin. The conformal mapping theorem \cite{41}(\S 2) associates to the interior of such a curve, a closed 2-cell in the sense of differential topology \cite{33}(\S 3.1), defining a `screw-top jar' map
\[
\tilde{\bl} :{\rm Annuli} \subset \haut (S^1) \to \Omega^2_e S^2
\]
% ? Dehn solenoid? 
from Segal's complexification of $\D S^1$ as the semigroup of annuli (recall $\D S^1 \supset \Sl_2(\R)$  \cite{29}) to the group completion of the monoid of homotopy cylinders. 
%(perhaps carrying a nonvanishing vector field)
\bigskip  

{\sc Proposal :} Birkhoff/Connes/Kreimer/Marcolli renormalization admits a generalization as a homotopy-theoretic pushout 
\[
\xymatrix{
B\Cz \ar[d]^\Lc \ar[r]^\rho & B (\Gamma  ? \sim \T) \ar[d] \\
BG_{\rm CK} \ar@{.>}[r]^{\bro_\Lc} \ar[ur] & \Omega S^2 }
\]
with 
\[
\xymatrix{
M\U^*(\om) \ar[d]^{Z_\shuffle} \ar[r] & MU^*(BG_{\rm CK}) \ar[d] \\
M\U^*(B\T) \ar[r] & MU^*B\Cz  .} 
\]
 
We think of $\bro$ as defining a complex line bundle $\rho_\Lc$ on $B\Cz$ and thus a characteristic  class $ch_{MU}(\Lc) \in MU^*(BG_{\rm CK})$. The proposed model lifts a \l{}ine bundle $\bro_\Lc$ to a complex line bundle $\rho_\Lc$ with cobordism Chern character
\[
ch(\rho_\Lc) =  Z_\shuffle \bro^*_\Lc \cz_{MU}(b) 
\]
renormalized by $Z_\shuffle$ as in Prop 2.3 above. \bigskip

{\bf appendix}  :  {\sc miscellaneous digressions} \bigskip

$\bvdash$ {\it i}) : Extracted from {\bf Loop Groups} : the Dirac Sea 

Let $\Z_+$ denote the set of nonnegative integers. A subset $S$ of $\Z$ is 
{\bf commensurable} with $\Z_+$ if the symmetric difference $S \; \Delta \; \Z_+$ is finite; in that case we can define 
\[
\# S := |S \; \Delta \; \Z_+| 
\]
to be the relative cardinality of $S$.  

This defines an integer-valued grading on the sets commensurable 
with $\Z_+$. Listing the elements of $S$ in increasing order defines a 
sequence $s_0 < s_1 < \cdots$ such that 
\[
s_k = k + \# \; S
\]
for all $k \geq k_0$, when $k_0$ is sufficiently large. If $\m$ is 
the least such $k_0$, then
\[
\pi(S) := (s_1 - s_0) + \cdots + (s_\m - s_{\m-1})
\]
is an ordered partition of $s_\m - s_0$.

Thus an element $s_0$, together with an ordered partition $\pi$, determines
$S$ uniquely; but the partition alone (ie without its ordering) is not
enough to determine $S$ (or even $\# \; S$). 

\newpage

$\bvdash$ {\it ii}) :  Example re ${\mathbb P}^1(\C) - \{0,1,\infty\}$ and \cite{37} :
\[
\Maps_*(S^1 \vee S^1,\Omega S^2 ) \simeq \Maps_* (S^1 \wedge (S^1 \vee S^1),S^2) \simeq \Maps_*(S^2,S^2)^{\times 2}
\]
so 
\[
\Pic_\wzw(BF^n) = \pi_0 \Maps_* (BF^n,\Omega S^2) \cong \Z^n 
\]
while 
\[
\Pic_\C(BF^n) = \pi_0 \Maps_*(BF^n,B\T) \cong H^2(\vee^n S^1,\Z) = 0 .
\] 
Similarly, 
\[
\Pic_\wzw(\U(2)) \cong \pi_0\Maps_*(S^1,\om) \times \pi_0\Maps_*(S^3,\om) \cong 
\]
\[
\cong \pi_2(S^2) \times \pi_4(S^2) \cong \Z \times \Z_2 
\]
(for $\U(2) \cong S^1 \times S^3$ as spaces, if not as groups), while $\Pic_\C(\U(2) = 0$. \bigskip

$\bvdash$ {\it iii}) : Fred $\Z"L$ Cohen's theorem, that the Pontrjagin ring $H_*({\rm SF}(n+1))$ is commutative \cite{11}, is commensurable with Borel's theorem about toruses in Lie groups. It suggests a homotopy-theoretic approach to manifold topology based on the series $\G(n) = \Sf(n+1)$ of $H$-spaces rather than the classical series $\{{\rm Gl}(n)\}$ (in terms of Poincar\'e duality cobordism and torsion characteristic classes for spherical fibrations \cite{36}(\S 8).

Cartan/Klein geometry \cite{50} considers moduli spaces of locally homogeneous structures built from relative Lie groups $K \subset G$.  A plausible homotopy-theoretic example is the assertion that 
\[
{\rm Gl} \supset {\mathbb O} \subset {\mathbb G}  = \Sf
\]
defines a Kervaire-Milnor-Sullivan-Lurie spectral moduli stack $[{\mathbb G} // {\mathbb O}]$ of differentiable structures on the sphere. Familiar low dimensional exceptional identifications suggest interpreting the inclusion
\[
\T \simeq {\rm Gl}(2) \simeq {\rm Cyl} \subset \Omega^2_e S^2 \simeq \Sf(2)
\]
as a specialization of that interpretation, in terms of block-like \l{}ine bundles \cite{26}. Indeed, it seems that an analogy
\[
{\{\mathbb O} : {\mathbb G}\} \Rightarrow \{ \T : \Omega^2_e S^2 \}
\]
makes sense as a directed Waldhausen square (as in \cite{16}(\S 3.2); recall that the mapping class group of a cylinder is infinite cyclic, generated by a Dehn \cite{51} twist). \bigskip

$\bvdash$ {\it iv}) : The Postnikov section
\[
\om  \to \om[0,2] \simeq |[B \T\sslash \Z]|  
\]
(with $\Z$ acting on $\T$ by $k,z \mapsto z^{(-1)^k}$) can be regarded as the classifying space of a topological groupoid. A generator 
\[
H^1(\Omega S^2,\Z) \cong \Z \ni w
\]
of its first cohomology group is something like the first Stiefel-Whitney class in differential topology, see below. 

Note that Maycock's Brauer-Wall group of graded continuous trace $C^*$ algebras \cite{6}
\[
MC(\Omega S^2) \cong H^1(\om,\Z_2) \tilde{\times} H^3(\om,\Z) \cong \Z_2 \times \Z 
\]
has a canonical (up to a Klein viergroup) generator defining a natural associated Morita class. \bigskip

$\bvdash {\it iv}^*$)  Galewski and Stern \cite{52} construct a fibration 
\[
\xymatrix{
{} & H(\Theta_0,4) \ar[r] & B\HL \ar[d] \\
\Omega S^2 \ar[ur] & \ar[l] X \ar@{.>}[ur] \ar[r] & B\Top }
\]
which classifies lifts of the structure bundle on a topological manifold $X$ to a simplicial homology block bundle \cite{54}. Here $H$ signifies an Eilenberg-MacLane space, $\Theta_0$ is the (not finitely generated \cite{53}) group of three-dimensional PL homology spheres (mod cobordism through four-dimensional acyclic manifolds with vanishing Rokhlin invariant). \bigskip

{\sc lemma} \cite{5}
\[
\Theta_0 \ni \theta \mapsto \gamma^2 b \cdot \theta \in H^4(\Omega S^2; \Theta_0)
\]
is an isomorphism. 

[Since the Pontryagin ring $H_*(\Omega S^2,\Z) \cong \Z[x]$.]\bigskip

{\sc proposition} {\it A class in ${\rm Pic}_{\bf wzw} (X)$ defines an action of the group $\Theta_0$ on such lifts (to a homology manifold structure on $X$), through the composition}
\[
X \to \Omega S^2 \to H(\Theta_0,4) .
\]

This suggests interpreting such a class as a kind of level structure on the topological manifold $X$. 
% [{\bf 18/10/023}]
\bigskip

$\bvdash$ {\it v}) :  A contact manifold $X$ of dimension $2k+1$ has a contact one-form $\alpha$ such that
$\alpha \wedge (d\alpha)^k \in \Omega^{2k+1}$ nowhere equals zero. Following Stong, an almost contact (oriented) manifold is defined by a lift of its orthogonal frame bundle to the split extension
 \[
 1 \to \U(k) \to C_{2k+1} \to \Z_2 \to 0
 \]
defined by complex conjugation of unitary matrices. The tangent bundle $\tau = \eta \oplus \xi$ of such a structured manifold can be presented as the sum of a line bundle $\xi$ and a $2n$-plane bundle $\eta$ such that $\xi \otimes \eta \cong \eta$, \cite{19}(appendix), \cite{43}(\S 6).

In particular when $k = 1$ we have a diagram
\[
\xymatrix{
B\T \ar[d] \ar[r] & B\T \ar[d] \\
\om [0,2] \cong  |[B\T \sslash \Z]|  \ar[d] \ar[r] &  |[B\T \sslash \Z_2]| = BC_3 \ar[d] \\
 B\Z \simeq \T \ar[r] & B\Z_2 \simeq \R P^\infty}
\]
of fibrations lifting classical $w_1$ by $w$.\bigskip

$\bvdash$ {\it vi}) : In 1963 Reinhart \cite{40} defined a notion of Lorentz cobordism 
\[
0 \to L^* \to \Omega^*_R \to \Omega^* \to 0 
\]
(of oriented manifolds together with a nonvanishing normal vector field on the interior of the cobordism between them), with a split exact sequence in which $L^{2k} = \Z$, while $L^{4k+1} = \Z_2$ and $L^{4k-1} = 0$. This is interpreted as suggesting that the Euler characteristic of an oriented manifold has an associated mod two invariant in dimensions congruent to one mod four defined by Kervaire's semi-characteristic 
\[
k(X) \equiv \sum \dim H^{2i}(X,\R) 
\]
(mod two), which is an invariant of suitable scissors congruence \cite{16}). 

Stong shows that a closed almost contact manifold bounds an almost complex manifold, and interprets $k(X)$ as the obstruction to $X = \partial Y?$ bounding a stably almost contact manifold $Y$. \bigskip

$\bvdash$ {\it vii}) :  Atiyah and Singer \cite{2} show that the symbol of a first-order Real (pseudo-) differential operator 
\[
D = \oplus (-1)^i [*,d] |_{\Omega^{2i}} : \Omega^{2*}(X) \to \Omega^{2*}(X) 
\]
(roughly, a skew-adjoint square root of $-\Delta$ on even forms) defines an element $\sigma(D) \in KR^{-1}(TX)$ with topological index equal to $k(X)$; where $TX$ has a Real structure defined by 
\[
T_x X \ni (x,v) \mapsto (x,-v) .
\]
This suggests the existence of a universal construction
\[
KR^{-1}(B{\rm O}(4k+1)^V) \to KR^{-1}(TX) 
\]
induced by the classifying map for the tangent bundle (with standard representation $V = \R \times \C^{2k}$) for comparison with
\[
KR^{-1}(B{\rm O}(4k+1)^V) \to KR^{-1}(BC_{2k+1}^V) ,
\]
which is beyond my limited competence in Real representation theory. In any case,\bigskip

$\bvdash$ {\it viii}) : The composition 
\[
X \to \om \to \om[0,2]  \simeq BC_3 \to BC_{2k+1}
\]
suggests interpreting a \l{}ine bundle on a $4k + 1$ manifold as defining a stably almost contact structure of some kind, with a characteristic real line bundle $\xi$ and a circle bundle $\eta$. Note that $\eta^2 \in K\rm{O}_2$ loops to an interesting element of $K{\rm O}^1(\om)$. \bigskip

$\bvdash$ {\it ix} :  In dimension $4k - 1$ the semicharacteristic and its relation to spin Lorentz cobordism is of cosmological interest: in particular when $d = 3$ \cite{14} this suggests the bilinear form 
\[
\Omega^1(X) \ni \alpha,\beta \mapsto \int_X \alpha \wedge d\beta 
\]
on one-forms, and methods from link calculus \cite{37}. \bigskip

$\bvdash$ {\it x} ) : If $\varepsilon_E$ is the Euler class of the canonical line bundle on the projectification ${\mathbb P}_\C E \to X$ of a complex vector bundle $E$ over $X$, then the norm\begin{footnote}{\ie with $t(T) = \sum_{i \geq 0} t_i T^{i+1}$ regarded as an $M\U^*(X)[[t_*]]$-module endomorphism of $M\U^*({\mathbb P}_\C E)[[t_*]]$}\end{footnote}
\[
c_{t_*}(E) = \sum c_I(E) \cdot t^I = \det t(\varepsilon_E) 
\]
defines the total Chern polynomial of $E$. Evaluation of Chern words $c_I$ of weight $|I|$ on the fundamental class of an $|I|$-dimensional manifold \cite{20}(\S 2.3) defines complete symmetric functions $h_I \in H_*(B\U)$, %? \cite{15}(\S 2, \S 4.9), 
which can be conveniently identified with their Hall duals $m_I \in H^*(B\U)$. \bigskip

$\bvdash$ {\it xi} : ) A linear functional $b : \g_\CK \to \C$ is, roughly, a physical `character', with an associated geometric character 
\[
\exp(-tb) : \G_\CK \to \C^\times
\]
(with a sign convenient for heat kernel issues). \bigskip

$\bvdash$ {\it xii}) :  There is a plausible analogy between noncommutative spans or correspondences based on $M\xi_*(X \times Y) \otimes \Q$ and the theory of distributions, via a Schwartz kernel/slant product lemma \cite{47}
\[
H_*(X \times Y,A') \to (H^*(Y,A) \otimes H^*(X,A))'
\]  
with $A = \qsymm$ as an algebra of smooth test functions, and its linear dual $A' = \nsymm$ as  generalized functions under convolution. \bigskip

$\bvdash$ {\it xiii}) : the Nambu-Goto area functional:

{\bf 1.1} Recall that an endomorphism $f : M \to M$ of a compact connected 
$n$-dimensional oriented closed smooth manifold sends its fundamental class 
$[M] \in H^n(M,\Z)$ to $f_*[M] = (\deg f) \cdot [M]$, for some $\deg f \in 
\Z$. The Kronecker dual class $[\dvol_M] \in H^n_\dR(M)$ in de Rham cohomology 
therefore satisfies
\[
\int_M f^*(\dvol_M) = \deg f \in \Z \;.
\]

Let us regard the unit 3-sphere $S^3 = B^3_+ \cup B^3_-$ as the union of 
two suitably oriented closed 3-balls overlapping on $S^2 = B^3_+ \cap B^3_-$,
and let $\phi_\pm \in \Maps(B^3_\pm,S^3)$ be smooth maps such that 
$\phi_+ = \phi_-$ in a collar neighborhood of their overlap. Let
\[
\phi := \phi_+ \cup \phi_- : B^3_+ \cup B^3_- \cong S^3 \to S^3
\]
be the map constructed by glueing $\phi_+$ and $\phi_-$ together on the 
boundaries of their domains, and let 
\[
\psi := \phi|_{S^2} \in \Omega^2 S^3
\]
be the restriction of the resulting map to that overlap. Then 
\[
\int_{S^3} \phi^* (\dvol_{S^3}) = \int_{B^3_+} \phi_+^* (\dvol_{B^3_+}) - 
\int_{B^3_-} \phi_-^* (\dvol_{B^3_-}) = \deg \phi \in \Z \;,
\]
so 
\[
\int_{B^3_+} \phi_+^* (\dvol_{B^3_+}) \equiv \int_{B^3_-} \phi_-^* 
(\dvol_{B^3_-})  \in \R {\rm mod}  \Z 
\]
is well-defined. In physics this is the classical Nambu-Goto functional
\[
\Omega^2 S^3 \ni \psi \mapsto \alpha(\psi) \in \R/\Z 
\]
which (roughly) measures  the area of the surface traced out by the motion of 
a closed string in three-space \cite{48}

% In lieu of unavailable graphics, consider the classical idea of the World as the upper hemisphere of a globe whose underside was the Ocean that the Sun sailed around to get to morning on the other side of the World. 
%
% The Wise knew that the World was in fact a spherical globe floating in the Empyrean, because they could see with their own eyes at the end of the month when the Sun and Full Moon were at the end of Day that the Moon was another spherical globe as is the World, and that it was illuminated by the Light of the Sun, its phases being Shadows... as the Sharp-Eyed saw as well with the Morning/Evening Star.
%
% What they did not know was what happened at the Equator; some thought it too fiery for Life but it was rumored that the Ocean pours over the Edge into the Abyss. Suppose instead that the lower hemisphere is entirely ocean, making the Sun's voyage more plausible. 
%
% In this language, the abstract geometric Equator, the edge of the Ocean bounding the southern hemisphere, is guarded by the MidGard Serpent who maintains the Horizon. As the Serpent writhes, the Horizon ripples across the Night Sky; its total angular fluctuation, suitably normalized, is the action in the Nambu-Goto model, see also Taiji. 

{\bf 1.2 Example} The Mercator projection 
\[
\mu : [-\pi,+\pi] \times S^2 \to S^3
\]
defines a family of loops which expand from a point and then collapse 
again to a point. More precisely,
\[
\lambda \mapsto \mu([-\pi,\lambda] \times S^2) := B_\lambda \subset S^3
\]
defines a family of closed 3-cells, with $\partial B_\lambda$ being the 
two-sphere at lattitude $\lambda \in (-\pi,+\pi)$. Giving the unit radius 
$n$-sphere its usual Riemannian metric, we have a one-to-one map
\[
(-\pi,+\pi) \ni\lambda \mapsto \pi^{-2} {\rm vol}(B_\lambda) \equiv 
\alpha (\partial B_\lambda) \in \R/\Z \;.
\]
The adjoint to $\mu$ then defines a composition
\[
\xymatrix{
\tilde{\mu} : S^1 \ar[r]^-{\tilde{\mu}} & \Omega^2 S^3 \ar[r]^-\alpha & S^1}
\] 
inducing isomorphisms
\[
\pi_1(S^1) \cong \Z \to \pi_1(\Omega^2 S^3) \cong \pi_3 S^3 \cong \Z \;.
\]  

{\bf 1.3} This construction extends more generally to define a function
\[
\alpha_G : \Omega^2 G \to \R/\Z 
\]
for any connected simple simply-connected Lie group, \eg $G = \SU(2) ( \cong S^3)$: Bott says that for such a group, $\pi_3(G)$ is infinite cyclic. A choice of generator defines a Kronecker dual cocycle 
\[
[\gamma_G : (u,v,w) \mapsto |[u,[v,w]|_K] \in \Lambda^3 \g^*
\]
(where $|-,-|_K$ is the Killing form) in the Chevalley-Eilenberg complex for 
the Lie algebra $\g$ of $G$ (which calculates the de Rham cohomology of 
$G$ as a space). The differential of $\phi_\pm : B^3_\pm \to G$ defines a 
map
\[
\phi^{-1}_\pm d\phi_\pm : B^3_\pm \to \g 
\]
identifying $\gamma_G(\phi^{-1}_\pm d\phi_\pm))$ with a three-form on 
$B^3_\pm$, and when $G$ is $\SU(2)$, this equals the pullback under $\phi_\pm$
of the usual volume form (defined by Haar measure) on $G = S^3$. By the 
argument above, its integral over $B^3_\pm$ is well-defined modulo $\Z$.

In fact we can replace $S^2 = \C P^1$ in this construction with any compact
Riemann surface $S = \partial W$, regarded as the boundary of a compact 
3-manifold, to define a similar functional 
\[
\alpha_{S,G} : \Maps(S,G) \to \R/\Z 
\]
on a space of base-pointed maps

{\bf 2.1} The Wess-Zumino-Witten model for string theory is defined on such 
infinite-dimensional spaces of maps in terms of Feynman measure 
\[
\exp(2\pi ik \Lc (\psi)) \D\psi = \exp (2 \pi ik \; \alpha_{S,G}
(\psi)) \cdot \exp (\pi i k \int_S |\psi^{-1} d \psi|^2_K) \: \cD\psi \;,
\]
where $k \in \Z$ is a parameter, something like an inverse to Planck's 
constant. The second term in the product (defined by an analog of 
kinetic energy) is a relatively well-behaved Gaussian measure, regarded as 
defining a kind of free field, while the first term is meaningful even though 
the Nambu-Goto term is defined only modulo $\Z$. The related Polyakov measure
agrees with $\alpha$ at critical points.

{\bf 2.2} Since $\pi_1 \Maps(S,G) \cong \pi_0 \Maps(\Sigma S,G)$, evaluating 
a map on homology defines a homomorphism
\[
\pi_1 \Maps(S,G) \to {\rm Hom}(H_3(\Sigma S),H_3(G)) \cong \Z \;.
\]
It is tempting to conjecture that this is an isomorphism; the 
example above implies this is true in the simplest case when $S$ is $S^2$ 
and $G$ is $S^3$. This is equivalent to the conjecture that the diagram 
\[
\xymatrix{
\Omega^2_{\rm reg}S^3 \ar@{.>}[d] \ar@{.>}[r] & \Z \times \Omega^2 S^3 \ar[d] 
\ar[ddr]^\alpha \\
{\rm Con}_\otimes(\C) \ar[drr]^-\Delta \ar[r] & \Omega^2 S^2 \\
{} & {} & \R/\Z }
\]
is homotopy-commutative: in which the upper square is a fiber product, defining
a space $\Omega^2_{\rm reg} S^3$ of `regularized' or `mollified' maps from 
$S^2$ to $S^3$ which go, under the Hopf map $S^3 \to S^2$, to the `nice' 
elements of $\Omega^2 S^2$ defined by Segal's charged particles. 

I think this says that the Nambu-Goto map is essentially homotopic to the
map defined by the discriminant, and that we can think of the configuration
space ${\rm Con}_\otimes$ as a nice space of states for the WZW model.\bigskip

\bibliographystyle{amsplain}

\end{document}